\documentclass[12pt]{article}

\usepackage{latexsym}
\usepackage{amsfonts}
\usepackage{amsmath}
\usepackage{amsthm}
\usepackage{amssymb}

\newtheorem{thm}{Theorem}[section]
\newtheorem{lem}{Lemma}[section]

\def\infint{\int_{-\infty}^\infty}

\def\convd{\stackrel{\cal D}{\rightarrow}}

\def\ex{{\rm E\,}}

\def\var{\mathop{\rm Var}\nolimits}

\begin{document}

\title {Asymptotic Normality of Nonparametric\\
 Kernel Type Deconvolution
Density Estimators: \\
crossing the Cauchy boundary}

\author {A.J. van Es \quad {\normalsize and}\quad H.-W. Uh\\[.3cm]
{\normalsize Korteweg-de Vries Institute for Mathematics,
 University of Amsterdam}\\
{\normalsize Plantage Muidergracht 24,
 1018 TV Amsterdam}\\
{\normalsize The Netherlands}}
\date{}
\maketitle

\begin{abstract}
We derive asymptotic normality of kernel type deconvolution density
estimators. In particular we consider deconvolution problems where the known
component of the convolution has a symmetric $\lambda$-stable
distribution with
$0<\lambda\leq 2$.
It turns out that the limit behavior changes if the exponent
parameter $\lambda$ passes the value one, the  case of Cauchy deconvolution.
\\[.5cm]
{\sl AMS classification:} primary 62G05; secondary 62E20\\[.1cm]
{\it Keywords:} deconvolution, kernel estimation,
asymptotic normality\\[.2cm]

\end{abstract}

\section{Introduction}

\setlength{\baselineskip}{0.7cm}

Let $X_1,\ldots, X_n$ be i.i.d. observations, where
$ X_i=Y_i+Z_i $ and
 $Y_i$ and $Z_i$ are independent random variables.
Assume that the unobservable $Y$'s have distribution function $F$ and density $f$,
and that the  $Z$'s have a known density $k$.
Note that $g$ equals $k*f$,
where $*$ denotes convolution.
The deconvolution problem is the problem of estimating the density $f$
 from the observations $X_i$ from the convolution density $g$.

A well known estimator of $f(x)$ is based on Fourier inversion
and kernel smoothing.
Let $w$ denote a {\em kernel function} and $h>0$ a {\em bandwidth}.
The kernel type estimator $f_{nh}(x)$ of the density $f$ at the point $x$ is defined as
\begin{equation}\label{estimator}
f_{nh}(x)=\frac{1}{2\pi} \int_{-\infty}^\infty
e^{-itx} \frac{\phi_w(ht)\phi_{emp}(t)}{ \phi_k(t)}\,dt
=
{1\over nh}\sum_{j=1}^n v_h\Big({{x-X_j}\over h}\Big),
\end{equation}
with
$$
v_h(u)={1\over 2\pi}\infint {{\phi_w(s)}\over\phi_k(s/h)}\
e^{-isu}ds.
$$
Here $\phi_{emp}$ denotes the empirical
characteristic function of the sample, i.e.
$
\phi_{emp}(t) = {1\over n}\sum_{j=1}^n e^{itX_j},
$
and $\phi_w$ and $\phi_k$ denote the characteristic functions of
$w$ and $k$ respectively. Note that, even though (\ref{estimator}) has the form of an
ordinary kernel density estimator, because of the dependence of
$v_h$ on the bandwidth $h$ it is different.
Kernel type estimators for the density $f$ and its distribution
function $F$ have been studied by many authors. Relatively recent papers
are Zhang (1990), Fan (1991a,b), Fan and Liu (1997),
Van Es and Kok (1998), Cator (2001), Van Es and Uh (2001),
and Delaigle and Gijbels (2002).
For an introduction see Wand and Jones (1995).
This paper covers a chapter in Uh (2003).

The expectation of the estimator (\ref{estimator}) has a familiar form.
We have, see for instance Stefanski and Carroll (1990),
\begin{equation}\label{expectation}
\ex f_{nh}(x)=\ex {1\over h}\,w\Big({{x-Y_j}\over
h}\Big).
\end{equation}
Indeed, this expectation is equal to the expectation of an ordinary kernel density estimator of
$f$ based on observations $Y_j$ from $f$.
Expansions of (\ref{expectation}) for $h\to 0$ are standard in kernel density
estimation theory and are hence omitted here. See for instance
Wand and Jones (1995).

Deconvolution problems are usually divided in  two groups, {\em ordinary
smooth deconvolution problems}, where the rate of decay to zero at infinity and minus infinity
of the characteristic function $\phi_k$ is algebraic, and {\em super smooth deconvolution
problems}, where it is essentially exponential. This rate of
decay, and hence the smoothness of the known density $k$, has a
tremendous influence on the variance of the estimator, see for instance Fan
(1991) or Cator (2001). By (\ref{expectation}) it is clear that the expectation is
not affected. The general picture is that with increasing smoothness of $k$ the
estimation problem becomes harder and the  the optimal rates become slower.

Our aim is to derive classical central limit type theorems for these kernel type
deconvolution estimators. For ordinary smooth deconvolution this
has first been achieved in Fan (1991) and extended in Fan and Liu (1997).
The limit behaviour in this case
is essentially equal to that of a kernel estimator of a higher order
derivative of a density.
In some specific deconvolution problems this is evident from
relatively simple inversion formulas, cf. Van Es and Kok (1998).
For instance, for  {\em generalized gamma deconvolution}  where $k$ is the density of
$\lambda_1E_1+\lambda_2E_2+\dots+\lambda_m E_m$,
with $\lambda_1>1,\dots,\lambda_m>0$ and $E_1,\dots,E_m$
independent standard exponential random variables we have
\begin{equation}\label{asnorsmooth}
\sqrt{nh^{2m+1}}(f_{nh}(x)-\ex f_{nh}(x))\convd N(0,\sigma^2),
\end{equation}
where $\sigma^2=(\lambda_1\dots\lambda_m)^2\int
w^{(m)}(v)^2dv\, g(x)$.
This result is typical for ordinary smooth deconvolution,
 in the sense of a rate of convergence that is algebraic in
$h$.

Asymptotic normality of $f_{nh}(x)$ in  super smooth deconvolution problems
has been derived by Zhang (1991), Fan (1991b) and Fan and Liu (1997).
Under suitable conditions their  theorems state
\begin{equation}
\sqrt{n}\, {{f_{nh}(x) - \ex f_{nh}(x)}\over s_n}\convd N(0,1),
\end{equation}
where $Z_{nj}={1\over h}\,v_h\Big({x-X_j\over h}\Big),\ j=1,\dots
n$
and either
$s_n^2={1\over n}\sum_{j=1}^n Z_{nj}^2$
or $s_n^2$ equals the sample variance of $Z_{n1},\dots,Z_{nn}$.
So the estimator is studentized in some respect. The asymptotic variance is not clear.
Van Es and Uh (2001)  derive a   central limit
type theorem like (\ref{asnorsmooth}) for super smooth deconvolution, where
the asymptotic variance is clear and the normalisation is deterministic.
Their result is given in Theorem \ref{supsmodens} below.

\bigskip

\noindent{\bf Condition W}

Let $\phi_w$ be real valued, symmetric and have support $[-1,1]$.  Let
$\phi_w(0)=1$, and let\\
$\phi_w(1-t)=At^\alpha + o(t^\alpha)$ as $t\downarrow 0$, for some constants $A$ and $\alpha\geq 0$.

\bigskip

\noindent{\bf Condition K}

Assume that $\phi_k$ has exponentially decreasing tails, i.e.
$
\phi_k(t)\sim C |t|^{\lambda_0} e^{-|t|^\lambda/\mu},
$
as $|t|\rightarrow \infty$, for some $\lambda>1, \mu>0, \lambda_0$, and
some real constant $C$.  Furthermore assume $\phi_k(t) \neq 0$ for all
$t$.

\bigskip

\noindent Note that Condition K excludes the Cauchy distribution and all
other distributions for which the tail of the characteristic function
decreases more slowly than $e^{-|t|}$.

\bigskip

\begin{thm}\label{supsmodens}
Assume Condition W, Condition K and   $\ex X^2<\infty$.
Then, as $n\to\infty$ and $h\to 0$,
\begin{equation}\label{esuh}
\frac{\sqrt{n}}{h^{\lambda(1+\alpha)+\lambda_0-1}
e^{\frac{1}{\mu h^\lambda}}}\,(f_{nh}(x)-\ex f_{nh}(x))\convd
N(0,\frac{A^2}{2\pi^2}(\mu/ \lambda)^{2+2\alpha}(\Gamma(\alpha+1))^2).
\end{equation}
\end{thm}

\bigskip

\noindent Surprisingly, the asymptotic variance is distribution
free, in the sense that it does not depend on $f$ or $x$.
The condition $\lambda>1$ was needed to ensure that remainder terms in the proof
of this theorem are asymptotically negligible.
Note also the condition that the second moment of the observations
is finite.

By studying deconvolution problems where the known distribution is
a symmetric stable distribution we will investigate the asymptotic behavior of the
kernel deconvolution estimators if the conditions of Theorem \ref{supsmodens}
are not satisfied.

\section{Deconvolution for symmetric stable densities}\label{sect:cauchy:general}
\setcounter{equation}{0}

Consider
deconvolution for   {\em symmetric stable densities} $k$
which have characteristic function
\begin{equation}\label{symstab}
\phi_k(t)=e^{-  |t|^\lambda/\mu},\quad \mu >0\ \mbox{and}\ 0<\lambda\leq
2.
\end{equation}
The condition $0<\lambda\leq 2$ is necessary to ensure that $k$ is a
density, cf.  Chung (1974). Hence the normal distribution is, in
some sense, extreme.
Note that for $\lambda$ equal to one $k$ is a Cauchy density.
The only symmetric stable distribution with finite
second moment is the normal distribution for which $\lambda$ equals two.
This implies that, unless $\lambda$ equals two, the second moment of
the observations will be infinite.
Hence the normal distribution is the only symmetric stable
distribution for which Theorem \ref{supsmodens} applies.
We will derive a limit behavior, similar to that described by
Theorem \ref{supsmodens}, for cases
where $\lambda$ is larger than one. Of even more interest are the cases
where $\lambda$ is equal to one, i.e. Cauchy deconvolution, or smaller than one. It turns out
that, while  crossing the Cauchy boundary, a different limit behavior appears.

\bigskip

For simplicity we
only consider the sinc kernel, defined by
\begin{equation}\label{sinc}
w(x) =\sin x/(\pi x)\quad  (\phi_w(t)=I_{[-1,1]}(t)).
\end{equation}
Results for a more general class of kernels $w$ are given in Uh
(2003).

First we give a heuristic derivation of the results which are rigorously
proved in Section \ref{proofs}.
Note that the estimator $f_{nh}$
can be rewritten as
\begin{equation}\label{densf}
f_{nh}(x)=
\frac{1}{\pi nh}\sum_{j=1}^n
\int_{0}^1 \cos\Big({s\Big(\frac{X_j-x}{h}}\Big)\Big)
e^{(s/h)^\lambda/\mu}ds.
\end{equation}

\medskip
\noindent Let $S$ denote a random variable, independent of the $X_j$, having
probability density $f_S$ given by
\begin{equation*}
f_S(s)=\frac{1}{c(h)}\,e^{(s/h)^\lambda/\mu}
\end{equation*}
on the interval $[0,1]$, where the normalization constant $c(h)$ is
given by
$
c(h)=\int_0^1e^{(s/h)^\lambda/\mu}ds.
$
Then (\ref{densf}) is equal to
\begin{equation}\label{ce}
f_{nh}(x)=
\frac{c(h)}{\pi nh}\sum_{j=1}^n
\ex (\cos((S/h)(X_j-x))|X_j)=
\frac{c(h)}{\pi nh}\sum_{j=1}^n
\ex_S \cos((S/h)(X_j-x)).
\end{equation}

It turns out that the asymptotics are greatly determined by the
asymptotics of the distribution of the random variable $S$.
Let us first consider its expectation and variance.
The following lemma gives expansions of the normalization constant, the expectation
of $S$ and the variance of $S$. Its proof is given in Section
\ref{proofslemmas}.

\begin{lem}\label{expansions}
For $0<\lambda\leq 2$ and $h\to 0$ we have
\begin{eqnarray}
c(h)&=&
\Big(\frac{\mu}{\lambda}h^\lambda +O(h^{2\lambda})\Big)e^{
(1/h)^\lambda/\mu},\label{chexp}\\
\ex S  &=&1-\frac{\mu}{\lambda}\,h^{\lambda}
+o(h^{\lambda}) ,\label{expexp}\\
\var S &=&\frac{\mu^2}{\lambda^2}\,h^{2\lambda}
+o(h^{2\lambda}).\label{varexp}
\end{eqnarray}
\end{lem}

\noindent These expansions suggest to normalize $S$ as follows. Write
\begin{equation}\label{En}
E_n=\frac{\lambda}{\mu}\frac{(S - 1)}{h^\lambda}.
\end{equation}
The density function, $f_{E_n}$ say, of $E_n$ is given by
$$
f_{E_n}(v)=
{\mu\over\lambda}\frac{h^\lambda e^{1/(\mu h^\lambda)}}{c(h)}\,
e^{{1/(\mu h^\lambda)}((1+(\mu/\lambda)h^{\lambda} v)^\lambda-1)}
{\bf I}_{[-(\lambda/\mu) h^{-\lambda},0]}(v).
$$
By Taylor expansion and Lemma \ref{expansions} it converges uniformly
on bounded intervals to $e^{v}{\bf I}_{(-\infty,0]}(v)$. 
This implies that $E_n$ 
converges in
distribution to $-E$ where $E$ denotes a standard exponential random
variable.

For the terms in (\ref{ce}) we have
\begin{align}
&\ex_S \cos((S/h)(X_j-x))=
\ex_S \cos(((1 + S - 1)/h)(X_j-x))\nonumber\\
&=
\ex_S \cos((1/h)(X_j-x))\cos(((S - 1)/h)(X_j-x))\nonumber\\
&\quad-
\ex_S \sin((1/h)(X_j-x))\sin(((S - 1)/h)(X_j-x))\nonumber\\
&=
\cos\Big({X_j-x\over h}\Big)
\ex_{E_n}\cos\Big({\mu \over  \lambda}\,h^{\lambda-1}
E_n(X_j-x)\Big)\label{term1}\\
&\quad-
 \sin\Big({X_j-x\over h}\Big)
\ex_{E_n}\sin\Big({\mu \over \lambda}\,h^{\lambda-1}
E_n(X_j-x)\Big).\label{term2}
\end{align}
It now becomes apparent that we may expect different asymptotics in
the cases $0<\lambda<1$, $\lambda=1$ and $1<\lambda\leq 2$. In these
cases the factor $h^{\lambda -1}$ in (\ref{term1}) and (\ref{term2})
diverges to infinity, equals one and vanishes.

\bigskip

The next three theorems establish asymptotic normality for
$1< \lambda\leq 2$, i.e.  for the
symmetric stable densities whose characteristic function
decreases more rapidly than the characteristic function of the
Cauchy distribution, for Cauchy deconvolution, and for
 $1/3< \lambda\leq 1$, i.e. for  the
symmetric stable densities whose characteristic function
decreases more slowly than the characteristic function of the
Cauchy distribution.

\begin{thm}\label{subcauchy}
Let $w$ be the sinc kernel (\ref{sinc}).
If $1<\lambda\leq 2$  then, as
$n\to\infty$ and $h\to 0$, we have
\begin{equation}\label{nor4}
\frac{\sqrt{n}}{h^{\lambda-1}
e^{\frac{1}{\mu h^\lambda}}}\,(f_{nh}(x)-\ex f_{nh}(x))\convd
N(0,\frac{1}{2\pi^2}(\mu/ \lambda)^2).
\end{equation}
\end{thm}

\begin{thm}\label{cauchy}
Let $w$ be the sinc kernel (\ref{sinc}).
If $\lambda$ equals one, i.e. Cauchy deconvolution,  then, as
$n\to\infty$ and $h\to 0$, we have
\begin{equation}\label{nor2}
\sqrt{n} e^{-\frac{1}{\mu h}}(f_{nh}(x)-\ex f_{nh}(x))\convd
N(0,\sigma^2),
\end{equation}
with
\begin{equation}\label{var2}
\sigma^2={1\over2\pi^2}
\int {\mu^2\over{1+\mu^2(u-x)^2}}\,g(u)du.
\end{equation}
\end{thm}

\begin{thm}\label{supcauchy}
Let $w$ be the sinc kernel (\ref{sinc}).
If $1/3<\lambda<1$  then, as
$n\to\infty, h\to 0$ and $nh\to\infty$, we have
\begin{equation}\label{nor3}
{\sqrt{n}\over h^{(\lambda-1)/2}e^{(1/h)^\lambda/\mu}}\,
(f_{nh}(x)-\ex f_{nh}(x))\convd
N(0,\sigma^2),
\end{equation}
with
\begin{equation}\label{var3}
\sigma^2={1\over 2\pi}{\mu\over\lambda}
\, g(x).
\end{equation}
\end{thm}

\medskip
The global picture we see from these three theorems is that for
$1/3<\lambda<1$, apart from the exponential rate of convergence, the
asymptotic variance resembles the asymptotic variance of a kernel
density estimator, in the sense that it depends on the value of $g$ at
the point $x$, as in (\ref{asnorsmooth}).  This is typical for smooth deconvolution
problems, though the rate of the variance
is exponential in $h$ and not algebraic.
For Cauchy deconvolution we see that the
asymptotic variance depends globally on $g$. For
$1<\lambda \leq 2$,
the estimator is asymptotically distribution free. It shares the asymptotics
of Theorem \ref{supsmodens}, even though the second moment of the
observations is infinite for $1<\lambda<2$. Concluding we see that
the restriction $\lambda>1$ in Theorem \ref{supsmodens} is essential and
that the finite second moment condition might not be.
Crossing the Cauchy boundary we get different asymptotics.

\section{Proofs}\label{proofs}
\setcounter{equation}{0}

\subsection{Basic lemma}

\begin{lem}\label{pairind}

Let $Y_h=(X-x)/h \bmod 2\pi$.
As $n\to\infty$ and $h\to 0$,
\begin{equation}\label{pairind1}
(X,Y_h)\convd (X,U),
\end{equation}
where $U$ is uniformly distributed on the interval $[0, 2\pi]$.
Moreover $X$ and $U$ are independent.

Assume  $0<\lambda<1$. Let $z$ be a bounded periodic function with period $2\pi$ and let
$\tilde w$ be a continuous and integrable function such that $\tilde w$ is
monotone in the tails. Then, as $n\to\infty$ and $h\to 0$,  we have
\begin{equation}\label{pairind2}
h^{ \lambda-1}\ex \Big(z\Big({X_j-x\over h}\Big)
\tilde w(h^{\lambda -1}(X_j-x))\Big)
\to
\frac{1}{2\pi}\,\int_0^{2\pi} z(u)du  \int \tilde w(u)du\,g(x).
\end{equation}
\end{lem}

\medskip
\noindent{\bf Proof}

\noindent
Note that the density $g=k*f$ of $X$ is continuous and bounded.  For $-\infty<u<\infty$,
$0\leq y<2\pi$ and $M<u$, we have by a Riemann sum approximation of the integral of $g$ over the interval
$[M,u]$,
\begin{eqnarray*}
\lefteqn{P(M< X\leq u, Y_h\leq y)=
\sum_{i:M<(2\pi i+y)h+x\leq u}\int_{2\pi ih+x}^{(2\pi i+y)h+x}g(t)\,dt+O(h)}\\
&=&
\sum_{i:M<(2\pi i+y)h+x\leq u} yh g(\xi_{i,h})+O(h)
=
\frac{y}{2\pi}\sum_{i:M<(2\pi i+y)h+x\leq u}2\pi h
g(\xi_{i,h})+O(h)\\
&=&
\frac{y}{2\pi} \int_{M}^u g(u)\,du+o(1)=\frac{y}{2\pi}\, (G(u)-G(M))+ o(1),
\end{eqnarray*}
where $\xi_{i,h}$ is a point on the interval
$[2i\pi h+x,2i\pi h+yh+x] \subset [2i\pi h+x,2(i+1)\pi h+x]$.

For arbitrary $\epsilon>0$ choose $M(\epsilon)$ such that $G(M(\epsilon))=P(X\leq
M(\epsilon))\leq \frac{1}{3}\, \epsilon$ and $n_0(\epsilon)$
such that for $n\geq n_0(\epsilon)$ we have
$|P(M(\epsilon)< X\leq u, Y_h\leq y) -
\frac{y}{2\pi}\, (G(u)-G(M(\epsilon)))|\leq \frac{1}{3}\, \epsilon$. Then, for
$n\geq n_0(\epsilon)$,
\begin{align*}
|P(X\leq& u, Y_h\leq y) -\frac{y}{2\pi}\, G(u)|\\
&
\leq
|P(X\leq u, Y_h\leq y)-P(M(\epsilon)< X\leq u, Y_h\leq y)|\\
&
\quad+
|P(M(\epsilon)< X\leq u, Y_h\leq y) -
\frac{y}{2\pi}\, (G(u)-G(M(\epsilon)))|\\
&
\quad+
|\frac{y}{2\pi}\, (G(u)-G(M(\epsilon)))-
\frac{y}{2\pi}\, G(u)|\\
&
\leq
P(X\leq M(\epsilon), Y_h\leq y)
+
\frac{1}{3}\, \epsilon
+
\frac{y}{2\pi}\, G(M(\epsilon))\leq \epsilon.
\end{align*}
Hence $\lim P(X\leq u, Y_h\leq y)=\frac{y}{2\pi}\, G(u)$, which
proves (\ref{pairind1}).

To prove (\ref{pairind2}) write
\begin{eqnarray*}
\lefteqn{h^{ \lambda-1}\ex  z \Big({X_j-x\over h}\Big)\tilde w(h^{\lambda
-1}(X_j-x))}\\
&=&
h^{ \lambda-1}\int z\Big({{u-x}\over h}\Big) \tilde w (h^{\lambda
-1}(u-x)) g(u)du\\
&=&
h^\lambda\int z(v)\tilde  w (h^\lambda v  ) g(x+hv)dv
\\
&=&
h^\lambda\sum_{i\in\mathbb Z}\int_0^{2\pi}
z(t+2\pi i) \tilde w (h^\lambda(t+2\pi i)) g(x+h(t+2\pi i) )dt\\
&=&
\frac{1}{2\pi}\int_0^{2\pi}z(t) \rho(t)dt,
\end{eqnarray*}
where
\begin{equation}
\rho(t)=2\pi h^\lambda \sum_{i\in\mathbb Z}
\tilde w(h^\lambda(t+2\pi i)) g(x+h(t+2\pi i)).
\end{equation}

With $\xi_i(t)=h^\lambda(t+2\pi i)$ for $i\in \mathbb Z$, we have
\begin{eqnarray}
\lefteqn{\rho(t)=2\pi h^\lambda \sum_{i\in\mathbb Z}
\tilde w(\xi_i(t)) g(x+h^{1-\lambda}\xi_i(t))}\nonumber\\
&=&
g(x)\, 2\pi h^\lambda \sum_{i\in\mathbb Z}
\tilde w(\xi_i(t))
+ 2\pi h^\lambda \sum_{i\in\mathbb Z}
\tilde w(\xi_i(t)) (g(x+h^{1-\lambda}\xi_i(t))-g(x)).
\label{ro}
\end{eqnarray}

Let $M>0$ be such that $|\tilde w|$ is increasing  on $(-\infty,-M]$ and
decreasing on $[M,\infty)$. Then, for $0\leq t\leq 2\pi$, we have
\begin{equation}
2\pi h^\lambda \sum_{i:\xi_i(t)\leq -M}
|\tilde w(\xi_i(t))|
\leq
2\pi h^\lambda \sum_{i:\xi_i(2\pi)\leq -M}
|\tilde w(\xi_i(2\pi))|+O(h^\lambda)
\leq
\int_{-\infty}^{-M}|\tilde w(u)|du + o(1)
\label{bound1}
\end{equation}
and
\begin{equation}
2\pi h^\lambda \sum_{i:\xi_i(t)\geq M}
|\tilde w(\xi_i(t))|
\leq
2\pi h^\lambda \sum_{i:\xi_i(0)\geq M}
|\tilde w(\xi_i(0))|+O(h^\lambda)
\leq
\int_M^\infty|\tilde w(u)|du + o(1).
\label{bound2}
\end{equation}
The convergence of the sum to the integral follows from the approximation
from below of $\tilde w$ by a step function and the dominated
convergence theorem. The $o(1)$ terms in (\ref{bound1}) and (\ref{bound2}) can be chosen
such that they do not depend on $t$.

Moreover, note that
\begin{equation}\label{conv}
2\pi h^\lambda \sum_{i:-M<\xi_i(t)< M}
|\tilde w(\xi_i(t))|\to \int_{-M}^{M} |\tilde w(u)|du,
\end{equation}
uniformly for $t$ in $[0,2\pi]$,
by the continuity of $\tilde w$ and Riemann sum approximation.

By the bounds (\ref{bound1}) and (\ref{bound2}), and the uniform
convergence in (\ref{conv}), one can show
\begin{equation}\label{conv2}
2\pi h^\lambda \sum_{i\in \mathbb Z}
|\tilde w(\xi_i(t))|\to \int_{-\infty}^\infty |\tilde w(u)|du,
\end{equation}
uniformly for $t$ in $[0,2\pi]$.
Similarly one can show
\begin{equation}\label{conv3}
2\pi h^\lambda \sum_{i\in \mathbb Z}
\tilde w(\xi_i(t))\to \int_{-\infty}^\infty \tilde w(u)du,
\end{equation}
uniformly for $t$ in $[0,2\pi]$, which implies that the first term
in (\ref{ro}) converges to $g(x)\, \int_{-\infty}^\infty \tilde
w(u)du$, uniformly for $t$ in $[0,2\pi]$.
Since $g$ is uniformly continuous we have
$g(x+h^{1-\lambda}\xi_i(t))-g(x)\to 0$, uniformly for $i$ and $t$ with
$-M<\xi_i(t)<M$. Using (\ref{bound1}), (\ref{bound2}) and
(\ref{conv2}) one can show that the second term in (\ref{ro})
vanishes, uniformly for $t$ in $[0,2\pi]$.
Hence
$\rho(t)\to g(x)\, \int_{-\infty}^\infty \tilde w(u)du$, uniformly
in $t$.

Finally we get
\begin{align}
h^{ \lambda-1}&\ex  z \Big({X_j-x\over h}\Big)\tilde w(h^{\lambda
-1}(X_j-x))=\frac{1}{2\pi}\int_0^{2\pi}z(t) \rho(t)dt\nonumber\\
&\to
\frac{1}{2\pi}\,\int_0^{2\pi}z(u) du\int_{-\infty}^\infty \tilde w(u)du
\, g(x),
\end{align}
which completes the proof of the Lemma.
\hfill$\Box$

\subsection{Proofs of Theorems \ref{subcauchy}, \ref{cauchy}
and \ref{supcauchy}}\label{prtheorem}

We can derive a bound on the error in substituting $-E$ for
$E_n$ in the terms (\ref{term1}) and (\ref{term2}).  The proof is
given in Section \ref{proofslemmas}.

\begin{lem}\label{approxs}
If $0<\lambda\leq 2$, as $n\to\infty$ and $h\to 0$, we have almost surely
\begin{equation*}
\Big|\ex_{E_n}\cos\Big({\mu\over\lambda}\,h^{\lambda -1}
E_n(X_j-x)\Big)
-
\ex_E\cos\Big(-{\mu\over\lambda}\,h^{\lambda -1}E(X_j-x)\Big)\Big|
=
O(h^\lambda)
\end{equation*}
and
\begin{equation*}
\Big|\ex_{E_n}\sin\Big({\mu\over\lambda}\,h^{\lambda -1}
E_n(X_j-x)\Big)
-
\ex_E\sin\Big(-{\mu\over\lambda}\,h^{\lambda -1}E(X_j-x)\Big)\Big|
=
O(h^\lambda),
\end{equation*}
where $E$ is a standard exponential random variable.
\end{lem}

We can now approximate $f_{nh}(x)-\ex f_{nh}(x)$.  We have
\begin{align}
&f_{nh}(x)-\ex f_{nh}(x)\nonumber\\
&=\frac{c(h)}{\pi nh}\sum_{j=1}^n
\Big(\ex_S \cos((S/h)(X_j-x))-\ex \ex_S
\cos((S/h)(X_j-x))\Big)\nonumber\\
&=
\frac{c(h)}{\pi nh}\sum_{j=1}^n
\Big(
\cos\Big({X_j-x\over h}\Big)\ex_{E_n}\cos\Big({\mu \over  \lambda}\,h^{\lambda
-1}
E_n(X_j-x)\Big)\nonumber\\
&\hspace{2.5cm}
-\sin\Big({X_j-x\over h}\Big)\ex_{E_n}\sin\Big({\mu \over \lambda}\,h^{\lambda
-1}
E_n(X_j-x)\Big)\Big)\nonumber\\
&\hspace{2.5cm}
-\ex\Big(\cos\Big({X_j-x\over h}\Big)\ex_{E_n}\cos\Big({\mu \over
\lambda}\,h^{\lambda -1}
E_n(X_j-x)\Big)\Big)\nonumber\\
&\hspace{2.5cm}
+\ex\Big(\sin\Big({X_j-x\over h}\Big)\ex_{E_n}\sin\Big({\mu \over
\lambda}\,h^{\lambda -1}
E_n(X_j-x)\Big)\Big)\Big)\nonumber\\
&=
\frac{c(h)}{\pi nh}\sum_{j=1}^n
\Big(
\cos\Big({X_j-x\over h}\Big)\ex_{E}\cos\Big(-{\mu \over  \lambda}\,h^{\lambda -1}
E(X_j-x)\Big)\nonumber\\
&\hspace{2.5cm}
-\sin\Big({X_j-x\over h}\Big)\ex_{E}\sin\Big(-{\mu \over \lambda}\,h^{\lambda -1}
E(X_j-x)\Big)\Big)\nonumber\\
&\hspace{2.5cm}
-\ex\Big(\cos\Big({X_j-x\over h}\Big)\ex_{E}\cos\Big(-{\mu \over
\lambda}\,h^{\lambda -1}
E(X_j-x)\Big)\Big)\nonumber\\
&\hspace{2.5cm}
+\ex\Big(\sin\Big({X_j-x\over h}\Big)\ex_{E}\sin\Big(-{\mu \over
\lambda}\,h^{\lambda -1}
E(X_j-x)\Big)\Big)\Big)\nonumber\\
&\hspace{5cm}
+
O_P\Big({1\over\sqrt{n}}\,h^{2\lambda-1}e^{(1/h)^\lambda/\mu}\Big).
\label{approxerror}
\end{align}
The order of the remainder term follows from the fact that it is
equal to
an average of independent terms, each of which is equal to the sum
of
$$
\frac{c(h)}{\pi h}\Big(
\cos\Big({X_j-x\over h}\Big)\ex_{E_n}\cos\Big({\mu \over  \lambda}\,h^{\lambda
-1}
E_n(X_j-x)\Big)
-\cos\Big({X_j-x\over h}\Big)\ex_{E}\cos\Big(- {\mu \over  \lambda}\,h^{\lambda
-1}
E(X_j-x)\Big)\Big)
$$
and
$$
-\frac{c(h)}{\pi h}\Big(
\sin\Big({X_j-x\over h}\Big)\ex_{E_n}\sin\Big({\mu \over  \lambda}\,h^{\lambda
-1}
E_n(X_j-x)\Big)
-\sin\Big({X_j-x\over h}\Big)\ex_{E}\sin\Big(- {\mu \over  \lambda}\,h^{\lambda
-1}
ES(X_j-x)\Big)\Big),
$$
minus their expectations.
The variances of these terms are of
order $c(h)^2O(h^{2\lambda})/(\pi^2 h^2)$
by Lemma \ref{approxs}, which  is of
order $O(h^{4\lambda-2}e^{2(1/h)^\lambda/\mu})$ by (\ref{chexp}).
Hence the variance of the average is of order
$O(\frac{1}{n}h^{4\lambda-2}e^{2(1/h)^\lambda/\mu})$,
which yields the order of the remainder term (\ref{approxerror}) by the Markov
inequality.

A straightforward computation yields
\begin{equation*}
\ex_E\cos\Big(-{\mu\over\lambda}\,h^{\lambda -1}E(X_j-x)\Big)
=
{\lambda^2\over{\lambda^2+\mu^2h^{2\lambda-2}(X_j-x)^2}}
=w_1(h^{\lambda -1}(X_j-x))
\end{equation*}
and
\begin{equation*}
\ex_E\sin\Big(-{\mu\over\lambda}\,h^{\lambda -1}E(X_j-x)\Big)
=
-\frac{\mu\lambda h^{\lambda-1}(X_j-x)}{\lambda^2+\mu^2h^{2\lambda-2}(X_j-x)^2}
=w_2(h^{\lambda -1}(X_j-x)),
\end{equation*}
where
$$
w_1(u)={\lambda^2\over{\lambda^2+\mu^2u^2}}\quad\mbox{and}\quad w_2(u)=- {\mu\lambda
u\over {\lambda^2+\mu^2u^2}}.
$$
Note that $w_1$ and $w_2$ are continuous bounded functions with
$w_1(u)^2 + w_2(u)^2=w_1(u)$. Note also that $w_2$ is not
integrable, so we can not apply Lemma \ref{pairind} directly for
$\tilde w=w_2$. However, since $w_2^\alpha$ is integrable for $\alpha>1$,
it turns out that we can circumvent this problem.

Define the random variables $V_{nj}$   as
\begin{align}
V_{nj}&=\cos\Big(\frac{X_j-x}{h}\Big)w_1(h^{\lambda -1}(X_j-x))
-\sin\Big(\frac{X_j-x}{h}\Big)w_2(h^{\lambda -1}(X_j-x))\nonumber\\
&=
\cos(Y_h)w_1(h^{\lambda -1}(X_j-x))
-\sin(Y_h)w_2(h^{\lambda -1}(X_j-x)).
\end{align}
Then
\begin{equation}\label{fasv}
f_{nh}(x)-\ex f_{nh}(x)=\frac{c(h)}{\pi h}\,\frac{1}{n}\sum_{j=1}^n (V_{n,j}-\ex
V_{n,j})+
O_P\Big({1\over\sqrt{n}}\,h^{2\lambda-1}e^{(1/h)^\lambda/\mu}\Big).
\end{equation}
To prove our three theorems we will   check the Lyapounov
condition for $\frac{1}{n}\sum_{j=1}^n (V_{n,j}-\ex V_{n,j})$ to be
asymptotically normal, i.e. for some $\delta >0$ we have to check
\begin{equation}\label{lyapounov}
\frac{\ex |V_{n,j}-\ex V_{n,j}|^{2+\delta}}{n^{\delta/2}(\var(V_{n,j}))^{1+\delta/2}}\to 0.
\end{equation}
We will check this condition for $\delta$ equal to two.
Note that by the inequality $|a+b|^p\leq
2^p(|a|^p+|b|^p), p\geq 0$ we have $\ex (V_{n,j}-\ex
V_{n,j})^4=\leq 4(\ex V_{n,j}^4+(\ex V_{n,j})^4)$.

\medskip

For $1<\lambda\leq 2$ the factor $h^{\lambda-1}$
vanishes. Hence, by Lemma \ref{pairind}, we have $(h^{\lambda-1}X,Y_h)\convd
(0,U)$. Since we are dealing with bounded continuous functions of
$(h^{\lambda-1}X,Y_h)$ we also have
\begin{align*}
\ex V_{nj}&=\ex (\cos(Y_h)w_1(h^{\lambda -1}(X_j-x))
-\sin(Y_h)w_2(h^{\lambda -1}(X_j-x)))\nonumber\\ &\to
\ex (\cos(U)w_1(0) -\sin(U)w_2(0))=0\\
\noalign{\noindent and}
\ex V_{nj}^4&=\ex (\cos(Y_h)w_1(h^{\lambda -1}(X_j-x))
-\sin(Y_h)w_2(h^{\lambda -1}(X_j-x)))^4\nonumber\\ &\to
\ex (\cos(U)w_1(0) -\sin(U)w_2(0))^4=\ex \cos(U)^4=\frac{3}{8}.
\end{align*}
The asymptotic variance is given by
\begin{equation*}
\var (V_{n,j})=\ex V_{n,j}^2 - (\ex V_{n,j})^2
\to
\ex (\cos(U)w_1(0) -\sin(U)w_2(0))^2=\ex \cos(U)^2=\frac{1}{2}.
\end{equation*}
Let us  check (\ref{lyapounov}) with
$\delta$ equal to two. Indeed we have
\begin{equation}
\frac{\ex |V_{n,j}-\ex V_{n,j}|^4}{n(\var(V_{n,j}))^2}
=\frac{O(1)}{n(\frac{1}{2}+o(1))^2}\to 0.
\end{equation}
This shows that $\frac{1}{n}\sum_{j=1}^n (V_{n,j}-\ex V_{n,j})$
and $\frac{c(h)}{\pi h}\,\frac{1}{n}\sum_{j=1}^n (V_{n,i}-\ex V_{n,i})$
are asymptotically normally distributed. The asymptotic variance
of $\frac{c(h)}{\pi h}\,\frac{1}{n}\sum_{i=1}^n (V_{n,j}-\ex
V_{n,j})$ is given by
\begin{equation}
\var\Big(\frac{c(h)}{\pi h}\,\frac{1}{n}\sum_{j=1}^n (V_{n,j}-\ex
V_{n,j})\Big)= \frac{1}{n}\,\frac{c(h)^2}{\pi^2
h^2}\,\var(V_{n,1})
\sim
\frac{1}{n}\,\frac{1}{2\pi^2}\frac{\mu^2}{\lambda^2}\,h^{2\lambda-2}\,e^{2
(1/h)^\lambda/\mu}.
\end{equation}
by Lemma \ref{expansions}.

\medskip

Now consider Cauchy deconvolution where $\lambda$ equals one.
By Lemma \ref{pairind}, since we are dealing with bounded continuous functions of
$(X,Y_h)$, we   have
\begin{align*}
\ex V_{nj}&=\ex (\cos(Y_h)w_1(X-x)
-\sin(Y_h)w_2(X-x))\nonumber\\ &\to
\ex (\cos(U)w_1(X-x) -\sin(U)w_2(X-x))=0\\
\noalign{and}
\ex V_{nj}^4&=\ex (\cos(Y_h)w_1(X-x)
-\sin(Y_h)w_2(X-x))^4\nonumber\\
&\to
\ex (\cos(U)w_1(X-x) -\sin(U)w_2(X-x))^4.
\end{align*}
The asymptotic variance of $V_{n,j}$ is given by
\begin{eqnarray*}
\lefteqn{\var (V_{n,j})=\ex V_{n,j}^2 - (\ex V_{n,j})^2}\\
&\to&
\ex (\cos(U)w_1(X-x) -\sin(U)w_2(X-x))^2\\
&=&\ex \cos(U)^2\ex w_1(X-x)^2+\ex \sin(U)^2\ex w_2(X-x)^2\\
&=&
\frac{1}{2}\,\ex w_1(X-x)= \frac{1}{2}\int {1\over{1+\mu^2(u-x)^2}}\,g(u)du.
\end{eqnarray*}
As above this shows that (\ref{lyapounov}) is satisfied for
$\delta$ equal to two.
Hence $\frac{1}{n}\sum_{j=1}^n (V_{n,j}-\ex V_{n,j})$
and $\frac{c(h)}{\pi h}\,\frac{1}{n}\sum_{j=1}^n (V_{n,i}-\ex V_{n,i})$
are asymptotically normally distributed. The asymptotic variance
of $\frac{c(h)}{\pi h}\,\frac{1}{n}\sum_{i=1}^n (V_{n,j}-\ex
V_{n,j})$ is given by
\begin{eqnarray}
\lefteqn{\var\Big(\frac{c(h)}{\pi h}\,\frac{1}{n}\sum_{j=1}^n (V_{n,j}-\ex
V_{n,j})\Big)}\nonumber\\
&=& \frac{1}{n}\,\frac{c(h)^2}{\pi^2
h^2}\,\var(V_{n,1})
\sim
\frac{1}{n}\,e^{2/(\mu h)}\,\frac{1}{2\pi^2}
\int {\mu^2\over{1+\mu^2(u-x)^2}}\,g(u)du.
\end{eqnarray}
by Lemma \ref{expansions}.

\medskip

Note that, if  $0<\lambda<1$, the factor   $h^{\lambda-1}$ diverges to infinity.
In this case we have
$$
\ex V_{n,j}^4 =
\sum_{l=0}^4 {4\choose l}\ex
\cos\Big(\frac{X_j-x}{h}\Big)^l(-1)^{4-l}
\sin\Big(\frac{X_j-x}{h}\Big)^{4-l}w_1(h^{\lambda -1}(X_j-x))^l
w_2(h^{\lambda -1}(X_j-x))^{4-l}.
$$
Since, for $l=0,1,\dots,4$, the functions $w_1^lw_2^{4-l}$ are integrable and monotone in the
tails, by Lemma \ref{pairind} we get
\begin{equation}\label{fourthpower}
\ex V_{n,j}^4 =O(h^{1-\lambda}).
\end{equation}
By a similar argument we have $\ex V_{n,j}^2=O(h^{1-\lambda})$, and hence  $(\ex V_{n,j})^4\leq (\ex
V_{n,j}^2)^2=O(h^{2-2\lambda})=O(h^{1-\lambda})$.

Next let us consider $\var (V_{n,j})=\ex V_{n,j}^2-(\ex
V_{n,j})^2$. Using the inequality above for $p=3/2$, the fact that $w_2^{3/2}$ is integrable,
 and Lemma \ref{pairind}, we get
$\ex |V_{n,j}|^{3/2}=O(h^{1-\lambda})$. By the Jensen
inequality we have
$|\ex V_{n,j}|^{3/2}\leq (\ex |V_{n,j}|)^{3/2}\leq  \ex
|V_{n,j}|^{3/2}$, and so
$(\ex V_{n,j})^2\leq ((\ex|V_{n,j}|^{3/2})^{2/3})^2=O(h^{(4/3)(1-\lambda}))$.
Moreover, by    (\ref{pairind2}) we get
\begin{equation*}
 \ex  (\cos\Big({X_j-x\over h}\Big)^2w_1(h^{\lambda -1}(X_j-x))^2
=
\frac{1}{2}\,h^{1-\lambda}\int w_1(u)^2dug(x)+o(h^{1-\lambda}).
\end{equation*}
Similarly we have
\begin{equation*}
\ex  (\sin\Big({X_j-x\over h}\Big)^2w_2(h^{\lambda -1}(X_j-x))^2
=
\frac{1}{2}\,h^{1-\lambda}\int w_2(u)^2dug(x)+o(h^{1-\lambda}).
\end{equation*}
and
\begin{eqnarray*}
\lefteqn{h^{\lambda -1}\ex  \cos\Big({X_j-x\over h}\Big)\sin\Big({X_j-x\over h}\Big)
w_1(h^{\lambda-1}(X_j-x))w_2(h^{\lambda-1}(X_j-x))}\\
&=&
\frac{1}{2} \int_0^{2\pi}\sin(u)\cos(u)du\int w_1(u)w_2(u) dug(x)+o(h^{1-\lambda})=
o(h^{1-\lambda}).
\end{eqnarray*}
Hence
\begin{eqnarray*}
\lefteqn{\ex V_{n,j}^2= \ex  \Big(
\cos\Big({X_j-x\over h}\Big)w_1(h^{\lambda -1}(X_j-x))
-\sin\Big({X_j-x\over h}\Big)w_2(h^{\lambda -1}(X_j-x))\Big)^2}\\
&=&
\frac{1}{2}\,h^{1-\lambda}\int
(w_1^2(u)+w_1^2(u))du g(x)+o(h^{1-\lambda})
=\frac{1}{2}\,h^{1-\lambda}\int
w_1(u)du g(x)+o(h^{1-\lambda})\\
&=&
\frac{1}{2}{\lambda\over\mu}\,\pi h^{1-\lambda}
\, g(x)+o(h^{1-\lambda}),
\end{eqnarray*}
and
\begin{equation}\label{varv}
\var(V_{n,j})=\frac{1}{2}{\lambda\over\mu}\,\pi h^{1-\lambda}
\, g(x)+o(h^{1-\lambda})+O(h^{(4/3)(1-\lambda}))=\frac{1}{2}{\lambda\over\mu}\,\pi h^{1-\lambda}
\, g(x)+o(h^{1-\lambda}).
\end{equation}
Finally we check (\ref{lyapounov}) with
$\delta$ equal to two. Indeed we have
\begin{equation}
\frac{\ex |V_{n,j}-\ex V_{n,j}|^4}{n\var(V_{n,j}))^2}
=\frac{O(h^{1-\lambda})}{n(\frac{1}{2}{\lambda\over\mu}\,\pi h^{1-\lambda}
\, g(x)+o(h^{1-\lambda}))^2}=O\Big(\frac{h^\lambda}{nh}\Big)\to 0.
\end{equation}
This shows that $\frac{1}{n}\sum_{j=1}^n (V_{n,j}-\ex V_{n,j})$
and $\frac{c(h)}{\pi h}\,\frac{1}{n}\sum_{j=1}^n (V_{n,i}-\ex V_{n,i})$
are asymptotically normally distributed. The asymptotic variance
of $\frac{c(h)}{\pi h}\,\frac{1}{n}\sum_{i=1}^n (V_{n,j}-\ex
V_{n,j})$ is given by
\begin{eqnarray*}
\lefteqn{\var\Big(\frac{c(h)}{\pi h}\,\frac{1}{n}\sum_{j=1}^n (V_{n,j}-\ex
V_{n,j})\Big)= \frac{1}{n}\,\frac{c(h)^2}{\pi^2
h^2}\,\var(V_{n,1})}\\
&\sim&
\frac{1}{n}\,\frac{\mu^2h^{2\lambda}}{\pi^2 \lambda^2h^2}\,e^{2
(1/h)^\lambda/\mu}\frac{1}{2}{\lambda\over\mu}\,\pi h^{1-\lambda}
\, g(x)
=\frac{1}{2\pi}\,\frac{\mu}{\lambda}\,\frac{1}{n}\, h^{\lambda-1} e^{2
(1/h)^\lambda/\mu}
\, g(x).
\end{eqnarray*}
by Lemma \ref{expansions} and (\ref{varv}).

\bigskip

It is easy to check that in all three cases the approximation error
(\ref{approxerror}) is
of smaller order than the asymptotic standard deviation in the theorems, provided
$\lambda>1/3$.
Hence this error is indeed negligible. \hfill$\Box$

\section{Proofs of the lemmas}\label{proofslemmas}
\setcounter{equation}{0}

\subsection{Proof of Lemma \ref{expansions}}

Note that, for $m=0,1,\cdots$, and any $0<\epsilon<1$,
$$
\int_0^\epsilon
s^me^{(s/h)^\lambda/\mu}ds=O(e^{(\epsilon/h)^\lambda/\mu})=o(h^{3\lambda}e^{(1/h
)^\lambda/\mu}).
$$
The exponent $3\lambda$ is fairly arbitrary but it suffices for our purposes.
By Lemma 3.5 of Van Es and Uh (2001) we have
$$
\int_\epsilon^1 s^{m-2\lambda}
e^{(s/h)^\lambda/\mu}ds
= {\mu\over\lambda}\,h^\lambda e^{(s/h)^\lambda/\mu} +o(h^\lambda
e^{(s/h)^\lambda/\mu}).
$$
By applying integration by parts twice we  get,
for $0<\epsilon<1$,
\begin{eqnarray*}
\lefteqn{\int_\epsilon^1 s^me^{(s/h)^\lambda/\mu}ds=
{\mu\over\lambda}\,h^\lambda\int_\epsilon^1
s^{m-\lambda+1}\Big({\lambda\over\mu}{1\over h^\lambda}s^{\lambda-1}
e^{(s/h)^\lambda/\mu}\Big)ds}\\
&=&
{\mu\over\lambda}\,
h^\lambda\Big(\Big[s^{m-\lambda+1}e^{(s/h)^\lambda/\mu}\Big]^1_\epsilon
-(m-\lambda+1)\int_\epsilon^1 s^{m-\lambda}
e^{(s/h)^\lambda/\mu}ds\Big)\\
&=&
{\mu\over\lambda}\,
h^\lambda\Big(\Big[s^{m-\lambda+1}e^{(s/h)^\lambda/\mu}\Big]^1_\epsilon
-(m-\lambda+1)\Big({\mu\over\lambda}\,
h^\lambda\Big(\Big[s^{m-2\lambda+1}e^{(s/h)^\lambda/\mu}\Big]^1_\epsilon\\
&&\quad\quad
-(m-2\lambda+1)\int_\epsilon^1 s^{m-2\lambda}
e^{(s/h)^\lambda/\mu}ds\Big)\Big)\Big)\\
&=&
{\mu\over\lambda}\,h^\lambda e^{(1/h)^\lambda/\mu}
-{\mu^2\over\lambda^2}(m-\lambda+1)h^{2\lambda}e^{(1/h)^\lambda/\mu}\\
&&\quad\quad
+{\mu^3\over\lambda^3}(m-\lambda+1)(m-2\lambda+1)h^{3\lambda}e^{(1/h)^\lambda/\mu}+o(h^{3\lambda}e^{(1/h)^\lambda/\mu}).
\end{eqnarray*}
Hence
\begin{eqnarray}
\lefteqn{\int_0^1 s^me^{(s/h)^\lambda/\mu}ds=
{\mu\over\lambda}h^\lambda e^{(1/h)^\lambda/\mu}
-{\mu^2\over\lambda^2}(m-\lambda+1)h^{2\lambda}e^{(1/h)^\lambda/\mu}}\nonumber\\
&&\quad\quad
+{\mu^3\over\lambda^3}(m-\lambda+1)(m-2\lambda+1)h^{3\lambda}
e^{(1/h)^\lambda/\mu}+o(h^{3\lambda}e^{(1/h)^\lambda/\mu}).
\label{mexpansion}
\end{eqnarray}
This expansion is used repeatedly in the remainder of the proof.

For $m=0$ we get
\begin{eqnarray}
\lefteqn{c(h)=\int_0^1 e^{(s/h)^\lambda/\mu}ds=
{\mu\over\lambda}h^\lambda e^{(1/h)^\lambda/\mu}
+{\mu^2\over\lambda^2}(\lambda-1)h^{2\lambda}e^{(1/h)^\lambda/\mu}}
\nonumber\\
&&\quad\quad
+{\mu^3\over\lambda^3}(1-\lambda)(1-2\lambda)h^{3\lambda}
e^{(1/h)^\lambda/\mu}+o(h^{3\lambda}e^{(1/h)^\lambda/\mu}),
\label{chexpansion}
\end{eqnarray}
which proves (\ref{chexp}).

Furthermore, using $(1+x)^{-1}=1-x+x^2+o(x^2)$
for $x\downarrow 0$, we have
\begin{eqnarray}
\lefteqn{\frac{1}{c(h)}={\lambda\over\mu}h^{-\lambda} e^{-(1/h)^\lambda/\mu}
\Big(1+{\mu\over\lambda}(\lambda-1)h^{\lambda}
+{\mu^2\over\lambda^2}(1-\lambda)(1-2\lambda)h^{2\lambda}
+o(h^{2\lambda}\Big)^{-1}}
\nonumber\\
&=&
{\lambda\over\mu}h^{-\lambda} e^{-(1/h)^\lambda/\mu}
\Big(1-{\mu\over\lambda}(\lambda-1)h^{\lambda}
-{\mu^2\over\lambda^2}(1-\lambda)(1-2\lambda)h^{2\lambda}
\nonumber\\
&&\quad\quad
+\Big({\mu\over\lambda}(\lambda-1)h^{\lambda}
+{\mu^2\over\lambda^2}(1-\lambda)(1-2\lambda)h^{2\lambda}\Big)^2+o(h^{2\lambda})\Big)
\nonumber\\
&=&
{\lambda\over\mu}h^{-\lambda} e^{-(1/h)^\lambda/\mu}\Big(
1-{\mu\over\lambda}(\lambda-1)h^{\lambda}
+{\mu^2\over\lambda^2}(\lambda-\lambda^2)h^{2\lambda}+o(h^{2\lambda})
\Big).\label{cinvexp}
\end{eqnarray}
Hence, by (\ref{mexpansion}) for $m=1$ and (\ref{cinvexp}),
\begin{eqnarray}
\lefteqn{\ex S = \frac{1}{c(h)}\int_0^1 s e^{(s/h)^\lambda/\mu}ds}\nonumber\\
&=&
\Big(1-{\mu\over\lambda}(\lambda-1)h^{\lambda}
+{\mu^2\over\lambda^2}(\lambda-\lambda^2)h^{2\lambda}+o(h^{2\lambda})\Big)
\Big(1
-{\mu\over\lambda}(2-\lambda)h^{\lambda}\nonumber\\
&&\quad\quad
+
{\mu^2\over\lambda^2}(2-\lambda)(2-2\lambda)h^{2\lambda}+o(h^{\lambda})\Big)\nonumber\\
&=&
1
-{\mu\over\lambda}h^{\lambda}+{\mu^2\over\lambda^2}(2-2\lambda)h^{2\lambda}+o(h^{2\lambda}),
\label{sexp}
\end{eqnarray}
which proves (\ref{expexp}).

Similarly, by (\ref{mexpansion}) for $m=2$ and (\ref{cinvexp}),
\begin{eqnarray}
\lefteqn{\ex S^2 = \frac{1}{c(h)}\int_0^1 s^2e^{(s/h)^\lambda/\mu}ds}\nonumber\\
&=&
\Big(1-{\mu\over\lambda}(\lambda-1)h^{\lambda}
+{\mu^2\over\lambda^2}(\lambda-\lambda^2)h^{2\lambda}
+o(h^{2\lambda})\Big)\nonumber\\
&&\quad\quad
\Big(1
-{\mu\over\lambda}(3-\lambda)h^{\lambda}
+{\mu^2\over\lambda^2}(3-\lambda)(3-2\lambda)h^{2\lambda}
+o(h^{2\lambda})\Big)\nonumber\\
&=&
1
-2
{\mu\over\lambda}h^{\lambda}+{\mu^2\over\lambda^2}(6-4\lambda)h^{2\lambda}+o(h^{2\lambda}).
\label{s2exp}
\end{eqnarray}
Finally, by (\ref{sexp}) and (\ref{s2exp}) we get
\begin{eqnarray*}
\lefteqn{\var S=\ex S^2-(\ex S)^2}\\
&=&
1-2{\mu\over\lambda}h^{\lambda}+{\mu^2\over\lambda^2}(6-4\lambda)h^{2\lambda}+o(h^{2\lambda})
- \Big(1
-{\mu\over\lambda}h^{\lambda}+{\mu^2\over\lambda^2}(2-2\lambda)h^{2\lambda}+o(h^{2\lambda})\Big)^2\\
&=&
1-{\mu^2\over\lambda^2}h^{2\lambda}+o(h^{2\lambda}),
\end{eqnarray*}
which proves (\ref{varexp}).
\hfill$\Box$

\subsection{Proof of Lemma \ref{approxs}}

Let $\epsilon_n=-h^{\lambda/2}/\log h$ denote a sequence of positive
(for $h<1$)
 numbers converging to zero.
Note that for $|t|$ small enough we have
\begin{equation*}
\Big|{1\over t}\,\Big((1+t)^\lambda -1\Big)- \lambda \Big|
\leq\lambda|\lambda-1||t|.
\end{equation*}
With $t=\mu h^\lambda v/\lambda$, for
$-\epsilon_n \lambda h^{-\lambda}/\mu\leq v\leq 0$, we have $-\epsilon_n\leq
t\leq 0$,
 and for $n$ large enough
\begin{equation*}
\Big|{\lambda\over\mu}\,{1\over h^\lambda
v}\,\Big(\Big(1+{\mu\over\lambda}\,
h^\lambda v\Big)^\lambda -1\Big)- \lambda \Big|
\leq\lambda|\lambda-1|{\mu\over\lambda}\,h^{\lambda}|v|,
\end{equation*}
and
\begin{align*}
\Big|{1\over\mu}&\,{1\over h^\lambda}\,\Big(\Big(1+{\mu\over\lambda}h^\lambda
v\Big)^\lambda -1\Big)- v \Big|
\leq |\lambda-1|{\mu\over\lambda}\,h^{\lambda}v^2\\
&\leq
|\lambda-1| {\lambda\over\mu}\, h^{-\lambda}\epsilon_n^2=o(1).
\end{align*}
This implies that for these values of $v$, and $n$ large enough,
\begin{eqnarray*}
\lefteqn{|f_{E_n}(v)e^{-v}-1|}\\
&=&
\Big|{\mu\over\lambda}\,{h^\lambda\over  c(h)}
e^{{1/(\mu h^\lambda)}}
\exp\Big(
{1\over\mu}\,{1\over h^\lambda}\,\Big(\Big(1+{\mu\over\lambda}h^\lambda
v\Big)^\lambda -1\Big)- v
\Big) -1\Big|\\
&\leq&
{\mu\over\lambda}\,{h^\lambda\over  c(h)}
e^{{1/(\mu h^\lambda)}}
\Big|
\exp\Big(
{1\over\mu}\,{1\over h^\lambda}\,\Big(\Big(1+{\mu\over\lambda}h^\lambda
v\Big)^\lambda -1\Big)- v
\Big) -1\Big|\\
&&+
|1-{\mu\over\lambda}\,{h^\lambda\over  c(h)}
e^{{1/(\mu h^\lambda)}}|\\
&\leq&
2\,{\mu\over\lambda}\,{h^\lambda\over  c(h)}
e^{{1/(\mu h^\lambda)}}
\Big|
{1\over\mu}\,{1\over h^\lambda}\,\Big(\Big(1+{\mu\over\lambda}h^\lambda
v\Big)^\lambda -1\Big)- v
\Big|+O(h^\lambda)\\
&\leq&
2\, {\mu\over\lambda}\,{h^\lambda\over  c(h)}
e^{{1/(\mu h^\lambda)}}
|\lambda-1|{\mu\over\lambda}h^{\lambda}v^2
+O(h^\lambda)\\
&\leq&
3{\mu\over\lambda}h^{\lambda}v^2
+O(h^\lambda),
\end{eqnarray*}
where the remainder terms do not depend on $v$.

Now note that
\begin{equation*}
\int_{-(\lambda/\mu) h^{-\lambda}}^{-\epsilon_n(\lambda/\mu) h^{-\lambda}}
|f_{E_n}(v)-e^v|dv\leq
\int_{-(\lambda/\mu) h^{-\lambda}}^{-\epsilon_n(\lambda/\mu) h^{-\lambda}}
f_{E_n}(v)dv
+e^{-\epsilon_n(\lambda/\mu)
h^{-\lambda}}
\end{equation*}
and, by Lemma \ref{expansions}, using $(1+s)^\lambda-1 \leq
(\lambda \wedge 1)s$ for $-1\leq s\leq 0$,
\begin{eqnarray*}
\lefteqn{\int_{-(\lambda/\mu) h^{-\lambda}}^{-\epsilon_n(\lambda/\mu) h^{-\lambda}}
f_{E_n}(v)dv=
{\mu\over\lambda}\frac{h^\lambda e^{1/(\mu h^\lambda)}}{c(h)}\,
\int_{-(\lambda/\mu) h^{-\lambda}}^{-\epsilon_n(\lambda/\mu) h^{-\lambda}}
e^{{1/(\mu h^\lambda)}((1+(\mu/\lambda)h^{\lambda}
v)^\lambda-1)}dv}\\
&=&
\frac{\lambda}{\mu}\,h^{-\lambda} (1+o(1))
\int_{-1}^{-\epsilon_n}
e^{{1/(\mu h^\lambda)}((1+s)^\lambda-1)}ds
\leq
\frac{\lambda}{\mu}\,h^{-\lambda} (1+o(1))
\int_{-1}^{-\epsilon_n}e^{{1/(\mu h^\lambda)} (\lambda\wedge
1)s}ds\\
&\leq&
\frac{\lambda}{\lambda\wedge 1}\,(1+o(1))e^{-\epsilon_n(\lambda/\mu)
(\lambda\wedge 1)h^{-\lambda}}=O(h^\lambda).
\end{eqnarray*}

Because the absolute value of the cosine is bounded
by one we get
\begin{align*}
\Big|\ex_{E_n}&\cos\Big({\mu\over\lambda}\,h^{\lambda -1}
E_n(X_j-x)\Big)
-
\ex_E\cos\Big(-{\mu\over\lambda}\,h^{\lambda -1}E(X_j-x)\Big)\Big|\\
&\leq
\int |f_{E_n}(v) -e^v I_{(-\infty,0]}(v)|dv\\
&=
\int_{-(\lambda/\mu) h^{-\lambda}}^0 |f_{E_n}(v) -e^v|dv
+  \int_{-\infty}^{-(\lambda/\mu) h^{-\lambda}}e^vdv\\
&=
\int_{-\epsilon_n(\lambda/\mu) h^{-\lambda}}^0 |f_{E_n}(v)e^{-v} -1|e^vdv
+\int_{-(\lambda/\mu) h^{-\lambda}}^{-\epsilon_n(\lambda/\mu) h^{-\lambda}}
|f_{E_n}(v) -e^v|dv
+  e^{-(\lambda/\mu) h^{-\lambda}}\\
&\leq
\int_{-\epsilon_n(\lambda/\mu) h^{-\lambda}}^0 \Big({3\mu\over\lambda}h^{\lambda}v^2
+O(h^\lambda)\Big)e^vdv+O(h^\lambda)=O(h^\lambda),
\end{align*}
which proves the first statement of the lemma.

The second statement can be proved similarly.
\hfil$\Box$

\bigskip

\bigskip

\noindent{\bf \Large Acknowledgment}
The research of the second author has been financed by the Netherlands
Organization for the Advancement of Scientific Research (NWO).

\bigskip

\noindent{\bf\Large References}

\begin{verse}

[1] E. Cator,
Deconvolution with arbitrary smooth kernels,
{\em  Statist. \& Probab. Lett.} {\bf 54}, (2001), 205--215.

[2] A. Delaigle  and I. Gijbels,
Comparison of data-driven bandwidth selection procedures in
deconvolution  kernel density estimation, to appear in Computational
Statistics and Data Analysis, (2002).

[3] K.L. Chung,
{\em A Course in Probability Theory},
Academic Press, London, 1974.

[4] A.J. van Es  and A.R. Kok,
Simple kernel estimators for certain nonparametric deconvolution problems,
{\em  Statistics \& Probability Letters} {\bf 39}, (1998), 151--160.

[5] A.J. van Es  and H.-W. Uh,
Asymptotic normality of kernel type deconvolution estimators,
Math. Preprint Series 01-26, Korte\-weg-de Vries Instituut voor Wiskunde,
Universiteit van Amsterdam, 2001.

[6] J. Fan,
On the optimal rates of convergence for nonparametric deconvolution problems,
{\em  Ann. Statist.} {\bf 19}, (1991a), 1257--1272.

[7] J. Fan,
Asymptotic normality for deconvolution kernel density estimators,
{\em  Sank\-hy\=a Ser. A} {\bf 53}, (1991b), 97--110.

[8] Y. Fan  and Y. Liu,
A note on asymptotic normality for deconvolution kernel density estimators,
{\em  Sankhy\=a Ser. A} {\bf 59}, (1997), 138--141.

[9] M.P. Wand  and M.C. Jones,
 {\em Kernel Smoothing},
 Chapman and Hall,
London, 1995.

[10] L. Stefanski and R.J. Carroll, Deconvoluting kernel density
estimators, {\em Statistics} {\bf 21}, (1990), 169--184.

[11] H.-W. Uh,
{\em Kernel Deconvolution},
PhD. Thesis, University of Amsterdam, 2003.

[12] C.H. Zhang,
Fourier methods for estimating mixing densities and distributions,
{\em Ann. Statist.} {\bf 18}, (1990), 806--831.

\end{verse}

\end{document}